\documentstyle[12pt]{amsart}
\begin{document}
  \title{Biharmonic maps into a Riemannian manifold of non-positive curvature}
  \title[Biharmonic maps into a non-positive curvature manifold]
  {Biharmonic maps into a Riemannian manifold of non-positive curvature}
\author{Nobumitsu Nakauchi}
  \address{Graduate School of Science and Engineering, \newline
  Yamaguchi University, 
  Yamaguchi, 753-8512, Japan}
  \email{nakauchi@@yamaguchi-u.ac.jp}
   \author{Hajime Urakawa}
  \address{Division of Mathematics, Graduate School of Information Sciences, Tohoku University, Aoba 6-3-09, Sendai, 980-8579, Japan}
  \curraddr{Institute for International Education, 
  Tohoku University, Kawauchi 41, Sendai 980-8576, Japan}
  \email{urakawa@@math.is.tohoku.ac.jp}
  \author{Sigmundur Gudmundsson} 
  \address{Mathematics, Faculty of Science, Lund University, Box 118, S-221 00 Lund, Sweden} 
  \email{Sigmundur.Gudmundsson@@math.lu.se}
  \thanks{
  Supported by the Grant-in-Aid for the Scientific Reserch, (C) No. 21540207, Japan Society for the Promotion of Science. 
  }
\maketitle
\begin{abstract}
  We study biharmonic maps between Riemannian manifolds with finite energy and finite bi-energy. We show that if the domain is complete and the target of non-positive curvature, then such a map is harmonic. 
  We then give applications to isometric immersions and horizontally conformal submersions. 
    \end{abstract}
\numberwithin{equation}{section}
\theoremstyle{plain}
\newtheorem{df}{Definition}[section]
\newtheorem{th}[df]{Theorem}
\newtheorem{prop}[df]{Proposition}
\newtheorem{lem}[df]{Lemma}
\newtheorem{cor}[df]{Corollary}
\newtheorem{rem}[df]{Remark}
\section{Introduction}
Harmonic maps play a central role in geometry, They are critical points of the energy functional 
$E(\varphi)=\frac12\int_M\vert d\varphi\vert^2\,v_g$ 
for smooth maps $\varphi$ of $(M,g)$ into $(N,h)$, and the Euler-Lagrange equation is that the tension filed 
$\tau(\varphi)$ vanishes. 
By extending notion of harmonic map, in 1983, J. Eells and L. Lemaire \cite{EL1} proposed a problem 
to consider the biharmonic maps which are, 
by definition, 
critical points of the bienergy functional
\begin{equation}
E_2(\varphi)=\frac12\int_M
\vert\tau(\varphi)\vert^2\,v_g.
\end{equation}
After G.Y. Jiang \cite{J} studied the first and second variation formulas of $E_2$, 
extensive studies in this area have been done
(for instance, see 
\cite{BFO}, \cite{CMP}, \cite{LO}, \cite{LO2},  \cite{MO1}, \cite{O1},  \cite{S1},
\cite{IIU2}, \cite{IIU},  \cite{II}, 
 etc.). Notice that harmonic maps are always biharmonic by definition. 
\par
For harmonic maps, it is well known that: 
\par
{\em If a domain manifold $(M,g)$ is complete and has non-negative Ricci curvature, and the sectional curvature of a target manifold $(N,h)$ is non-positive, then 
every energy finite harmonic 
map is a constant map} (cf. \cite{SY}). 
\par
Therefore, it is a natural question to consider 
biharmonic maps into a Riemannian manifold of non-positive curvature. 
In this connection, Baird, Fardoun and Ouakkas (cf. \cite{BFO}) showed that:  
\par
{\em If a non-compact Riemannian manifold 
$(M,g)$ is complete and has non-negative Ricci curvature and $(N,h)$ has non-positive sectional curvature, then every bienergy finite biharmonic map of $(M,g)$ into 
$(N,h)$ is harmonic}.  
\vskip0.3cm\par
In this paper, we will show that 
\begin{th} $($cf. {\bf Theorem 2.1}$)$
Under only the assumptions of completeness of 
$(M,g)$ and non-positivity of curvature of $(N,h)$,  \par
$(1)$ every biharmonic map 
$\varphi:\,(M,g)\rightarrow (N,h)$ with 
finite energy and 
finite bienergy 
must be harmonic. 
\par
$(2)$ In the case ${\rm Vol}(M,g)=\infty$, under the same assumtion,  
every biharmonic map 
$\varphi:\,(M,g)\rightarrow (N,h)$ with finite bienergy 
is harmonic. 
\end{th}
\vskip0.3cm\par
We do not need any assumption on the Ricci 
curvature of $(M,g)$ in Theorem 1.1. 
Since $(M,g)$ is a non-compact complete Riemannian manifold whose Ricci curvature 
is non-negative, then ${\rm Vol}(M,g)=\infty$ 
(cf. Theorem 7, p. 667, \cite{Y}).   
Thus, Theorem 1.1, $(2)$ recovers the result 
of Baird, Fardoun and Ouakkas. 
Furthermore, Theorem 1.1 is sharp 
because one can not 
weaken the assumptions
because 
the generalized Chen's conjecture does not hold 
if $(M,g)$ is not complete (cf. recall the counter examples of Ou and Tang \cite{OT}).  
The both assumptions of finitenesses
of the energy and bienergy are also necessary. 
Indeed, 
there exists a biharmonic map $\varphi$ which is not harmonic,  
but energy and bienergy are infinite. For example, 
$f(x)=r(x)^2=\sum_{i=1}^m
(x_i{})^2, x=(x_1,\cdots,x_m)\in {\mathbb R}^m$ 
is biharmonic, but not harmonic,
and have 
infinite energy and bienergy. 
\par
As the first bi-product of our method, we obtain 
(cf. \cite{NU1}, \cite{NU2}) 
\begin{th} $($cf. {\bf Theorem 3.1}$)$
Assume that $(M,g)$ is a complete 
Riemannian manifold, and let 
$\varphi:\,(M,g)\rightarrow (N,h)$ 
is an isometric immersion, and the sectional curvature of $(N,h)$ is non-positive. 
If $\varphi:\,(M,g)\rightarrow (N,h)$ is biharmonic 
and $\int_M\vert\xi\vert^2\,v_g<\infty$, then 
it is minimal. 
Here, $\xi$ is the mean curvature normal vector field of the isometric immersion $\varphi$. 
\end{th}
\vskip0.3cm\par
Theorem 1.2 (cf. Theorem 3.1$)$ 
gives an affirmative answer to the generalized B.Y. Chen's conjecture (cf. \cite{CMP})
under natural conditions. 
\par
For the second bi-product, we can apply Theorem 1.1 to a horizontally conformal submersion (cf. \cite{BE},\cite{BW}). Then, we have
\begin{th} $($cf. {\bf Corollary 3.4}$)$ 
Let $(M^m,g)$ be a non-compact complete Riemannian manifold $(m>2)$, and $(N^2,h)$, a Riemannian surface with non-positive curvature.  
Let  $\lambda$ be a positive function 
on $M$ belonging to 
 $C^{\infty}(M)\cap L^2(M)$, 
 and 
 $\varphi:\, (M,g)\rightarrow (N^2,h)$, 
 a horizontally conformal submersion with a 
 dilation $\lambda$.  
If $\varphi$ is biharmonic and $\lambda\,\vert\hat{\mathbf H}\vert_g\in L^2(M)$, then 
$\varphi$ is a harmonic morphism. 
Here, $\hat{\mathbf H}$ is trace of the second fundamental form of each fiber of $\varphi$. 
\end{th} 
\vskip0.6cm\par
{\bf Acknowledgement.} \quad 
The second author would like to express his sincere gratitude to Professor Sigmundur Gudmundsson for his hospitality and very intensive and helpful discussions 
during 
for the second author staying at Lund University, 
at 2012 May. 
Addition of part of Section Three 
to the original version 
of the first and second authors was 
based by this joint work with him during this period. 
\vskip0.6cm\par
\section{Preliminaries and statement of main theorem}
In this section, we prepare materials for the first and second variational formulas for the bienergy functional and biharmonic maps. 
Let us recall the definition of a harmonic map $\varphi:\,(M,g)\rightarrow (N,h)$, of a compact Riemannian manifold $(M,g)$ into another Riemannian manifold $(N,h)$, 
which is an extremal 
of the {\em energy functional} defined by 
$$
E(\varphi)=\int_Me(\varphi)\,v_g, 
$$
where $e(\varphi):=\frac12\vert d\varphi\vert^2$ is called the energy density 
of $\varphi$.  
That is, for any variation $\{\varphi_t\}$ of $\varphi$ with 
$\varphi_0=\varphi$, 
\begin{equation}
\frac{d}{dt}\bigg\vert_{t=0}E(\varphi_t)=-\int_Mh(\tau(\varphi),V)v_g=0,
\end{equation}
where $V\in \Gamma(\varphi^{-1}TN)$ is a variation vector field along $\varphi$ which is given by 
$V(x)=\frac{d}{dt}\vert_{t=0}\varphi_t(x)\in T_{\varphi(x)}N$, 
$(x\in M)$, 
and  the {\em tension field} is given by 
$\tau(\varphi)
=\sum_{i=1}^mB(\varphi)(e_i,e_i)\in \Gamma(\varphi^{-1}TN)$, 
where 
$\{e_i\}_{i=1}^m$ is a locally defined frame field on $(M,g)$, 
and $B(\varphi)$ is the second fundamental form of $\varphi$ 
defined by 
\begin{align}
B(\varphi)(X,Y)&=(\widetilde{\nabla}d\varphi)(X,Y)\nonumber\\
&=(\widetilde{\nabla}_Xd\varphi)(Y)\nonumber\\
&=\overline{\nabla}_X(d\varphi(Y))-d\varphi(\nabla_XY),
\end{align}
for all vector fields $X, Y\in {\frak X}(M)$. 
Here, 
$\nabla$, and
$\nabla^N$, 
 are connections on $TM$, $TN$  of $(M,g)$, $(N,h)$, respectively, and 
$\overline{\nabla}$, and $\widetilde{\nabla}$ are the induced ones on $\varphi^{-1}TN$, and $T^{\ast}M\otimes \varphi^{-1}TN$, respectively. By (2.1), $\varphi$ is harmonic if and only if $\tau(\varphi)=0$. 
\par
The second variation formula is given as follows. Assume that 
$\varphi$ is harmonic. 
Then, 
\begin{equation}
\frac{d^2}{dt^2}\bigg\vert_{t=0}E(\varphi_t)
=\int_Mh(J(V),V)v_g, 
\end{equation}
where 
$J$ is an elliptic differential operator, called the
{\em Jacobi operator}  acting on 
$\Gamma(\varphi^{-1}TN)$ given by 
\begin{equation}
J(V)=\overline{\Delta}V-{\mathcal R}(V),
\end{equation}
where 
$\overline{\Delta}V=\overline{\nabla}^{\ast}\overline{\nabla}V
=-\sum_{i=1}^m\{
\overline{\nabla}_{e_i}\overline{\nabla}_{e_i}V-\overline{\nabla}_{\nabla_{e_i}e_i}V
\}$ 
is the {\em rough Laplacian} and 
${\mathcal R}$ is a linear operator on $\Gamma(\varphi^{-1}TN)$
given by 
${\mathcal R}(V)=
\sum_{i=1}^mR^N(V,d\varphi(e_i))d\varphi(e_i)$,
and $R^N$ is the curvature tensor of $(N,h)$ given by 
$R^N(U,V)=\nabla^N{}_U\nabla^N{}_V-\nabla^N{}_V\nabla^N{}_U-\nabla^N{}_{[U,V]}$ for $U,\,V\in {\frak X}(N)$.   
\par
J. Eells and L. Lemaire \cite{EL1} proposed polyharmonic ($k$-harmonic) maps and 
Jiang \cite{J} studied the first and second variation formulas of biharmonic maps. Let us consider the {\em bienergy functional} 
defined by 
\begin{equation}
E_2(\varphi)=\frac12\int_M\vert\tau(\varphi)\vert ^2v_g, 
\end{equation}
where 
$\vert V\vert^2=h(V,V)$, $V\in \Gamma(\varphi^{-1}TN)$.  
\par
Then, the first variation formula of the bienergy functional 
is given $($the first variation formula$)$ by
\begin{equation}
\frac{d}{dt}\bigg\vert_{t=0}E_2(\varphi_t)
=-\int_Mh(\tau_2(\varphi),V)v_g.
\end{equation}
Here, 
\begin{equation}
\tau_2(\varphi)
:=J(\tau(\varphi))=\overline{\Delta}(\tau(\varphi))-{\mathcal R}(\tau(\varphi)),
\end{equation}
which is called the {\em bitension field} of $\varphi$, and 
$J$ is given in $(2.4)$.  
\par
A smooth map $\varphi$ of $(M,g)$ into $(N,h)$ is said to be 
{\em biharmonic} if 
$\tau_2(\varphi)=0$. 
\vskip0.6cm\par
Then, we can state our main theorem. 
\begin{th} 
Assume that $(M,g)$ is complete and the sectional curvature of $(N,h)$ 
is non-positive. \par
$(1)$ Every biharmonic map 
$\varphi:\,(M,g)\rightarrow (N,h)$ with 
finite energy $E(\varphi)<\infty$ and 
finite bienergy $E_2(\varphi)<\infty$, 
is harmonic. \par
$(2)$ In the case ${\rm Vol}(M,g)=\infty$,  
every biharmonic map 
$\varphi:\,(M,g)\rightarrow (N,h)$ with finite bienergy $E_2(\varphi)<\infty$,
is harmonic. 
\end{th}
\vskip0.6cm\par
\section{Proof of Main Theorem and two applications}
In this section we will give a proof of 
Theorem 2.1 which consists of 
four steps. 
\par
({\it The first step}) \quad 
For a fixed point $x_0\in M$, and 
for every 
$0<r<\infty$, 
we first take a cut-off  $C^{\infty}$ function $\eta$ on $M$ 
(for instance, see \cite{K}) satisfying that 
\begin{equation}
\left\{
\begin{aligned}
0\leq &\eta(x)\leq 1\quad (x\in M),\\
\eta(x)&=1\qquad\quad (x\in B_r(x_0)),\\
\eta(x)&=0\qquad\quad (x\not\in B_{2r}(x_0)),\\
\vert\nabla\eta\vert&\leq\frac{2}{r}
\qquad\,\,\, (x\in M).
\end{aligned}
\right.
\end{equation}
\par
\vskip0.6cm\par
For a biharmonic map
$\varphi:\,(M,g)\rightarrow (N,h)$, 
the bitension field is given as 
\begin{equation}
\tau_2(\varphi)=
\overline{\Delta}(\tau(\varphi))
-\sum_{i=1}^m
R^N(\tau(\varphi),d\varphi(e_i))d\varphi(e_i)=0,
\end{equation}
so we have 
\begin{align}
\int_M\langle\overline{\Delta}(\tau(\varphi)),\eta^2\,\tau(\varphi)\rangle\,v_g
&=\int_M
\eta^2\sum_{i=1}^m
\langle R^N(\tau(\varphi),d\varphi(e_i))d\varphi(e_i),\tau(\varphi)\rangle\,v_g\nonumber\\
&\leq 0, 
\end{align}
since the sectional curvature of $(N,h)$ is non-positive.  
\vskip0.6cm\par
({\it The second step})
Therefore, by (3.3) and noticing that $\overline{\Delta}=\overline{\nabla}^{\ast}\,\overline{\nabla}$, 
we obtain 
\begin{align}
0&\geq\int_M\langle\overline{\Delta}(\tau(\varphi)), 
\eta^2\,\tau(\varphi)\rangle\,v_g\nonumber\\
&=\int_M\langle\overline{\nabla}\tau(\varphi),
\overline{\nabla}(\eta^2\,\tau(\varphi))\rangle\,v_g\nonumber\\
&=\int_M\sum_{i=1}^m\langle \overline{\nabla}_{e_i}\tau(\varphi),\overline{\nabla}_{e_i}(\eta^2\,\tau(\varphi))\rangle\,v_g\nonumber\\
&=\int_M
\sum_{i=1}^m\bigg\{
\eta^2\,\langle\overline{\nabla}_{e_i}\tau(\varphi),
\overline{\nabla}_{e_i}\tau(\varphi)\rangle+
e_i(\eta^2)\,\langle\overline{\nabla}_{e_i}\tau(\varphi),\tau(\varphi)\rangle\bigg\}\,v_g
\nonumber\\
&=\int_M\eta^2\sum_{i=1}^m\bigg\vert\overline{\nabla}_{e_i}\tau(\varphi)\bigg\vert^2\,v_g
+2\int_M\sum_{i=1}^m
\langle\eta\overline{\nabla}_{e_i}\tau(\varphi), 
e_i(\eta)\,\tau(\varphi)\rangle\,v_g, 
\end{align}
where we used $e_i(\eta^2)=2\eta\,e_i(\eta)$ at the last equality.  
By moving the second term in the last equality of (3.4) to the left hand side,  
we have 
\begin{align}
\int_M
\eta^2\sum_{i=1}^m\vert\overline{\nabla}_{e_i}\tau(\varphi)\vert^2\,v_g
&\leq -2\int_M\sum_{i=1}^m
\langle
\eta\,\overline{\nabla}_{e_i}\tau(\varphi),e_i(\eta)\,\tau(\varphi)\rangle\,v_g\nonumber\\
&=-2\int_M\sum_{i=1}^m\langle V_i,W_i\rangle\,v_g, 
\end{align}
where we put 
$V_i:=\eta\,\overline{\nabla}_{e_i}\tau(\varphi)$, and $W_i:=e_i(\eta)\,\tau(\varphi)$ ($i=1\,\cdots,m$).  
 \par
 Now let recall the following 
 Cauchy-Schwartz inequality: 
 \begin{equation}
 \pm 2\,\langle V_i,W_i\rangle
 \leq \epsilon \vert V_i\vert ^2+\frac{1}{\epsilon}\vert W_i\vert^2
 \end{equation} 
 for all positive $\epsilon>0$ because of the inequality 
 $
 0\leq \vert \sqrt{\epsilon}\,V_i\pm\frac{1}{\sqrt{\epsilon}}\,W_i\vert^2.  
 $
 Therefore, for (3.5), we obtain 
\begin{align}
-2\int_M\sum_{i=1}^m\langle V_i,W_i\rangle\,v_g
\leq
\epsilon\int_M
\sum_{i=1}^m\vert V_i\vert^2\,v_g
+\frac{1}{\epsilon}\int_M\sum_{i=1}^m\vert W_i\vert^2\,v_g. 
\end{align} 
If we put $\epsilon=\frac12$, we obtain, 
by (3.5) and (3.7), 
\begin{align}
\int_M
\eta^2\sum_{i=1}^m\vert\overline{\nabla}_{e_i}\tau(\varphi)\vert^2\,v_g
&\leq \frac12 
\int_M\sum_{i=1}^m\eta^2\,\vert\overline{\nabla}_{e_i}\tau(\varphi)\vert^2\,v_g
\nonumber\\
&\qquad+2\int_M\sum_{i=1}^me_i(\eta)^2\,\vert\tau(\varphi)\vert^2\,v_g.
\end{align}
Thus, by (3.8) and (3.1), we obtain 
\begin{align}
\int_M\eta^2\sum_{i=1}^m\vert\overline{\nabla}_{e_i}\tau(\varphi)\vert^2\,v_g
&\leq 
4\int_M\vert\nabla\eta\vert^2\,\vert\tau(\varphi)\vert^2\,v_g\nonumber\\
&\leq\frac{16}{r^2}\int_M\vert\tau(\varphi)\vert^2\,v_g.
\end{align}
\vskip0.6cm\par
({\it The third step})\quad 
Since $(M,g)$ is complete and non-compact, we can tend $r$ to infinity. 
By the assumption $E_2(\varphi)=\frac12\int_M \vert\tau(\varphi)\vert^2\,v_g<\infty$, 
the right hand side goes to zero. And also, if $r\rightarrow\infty$, the left hand side of (3.9) goes to 
$
\int_M\sum_{i=1}^m\vert\overline{\nabla}_{e_i}\tau(\varphi)\vert^2\,v_g$
since
$\eta=1$ on $B_r(x_0)$. 
Thus, we obtain 
\begin{align}
\int_M\sum_{i=1}^m\vert\overline{\nabla}_{e_i}\tau(\varphi)\vert^2\,v_g=0. 
\end{align}
Therefore, we obtain, for every vector field $X$ in $M$,  
\begin{align}
\overline{\nabla}_X\tau(\varphi)=0. 
\end{align}
Then, we have, in particular, 
$\vert \tau(\varphi)\vert$ 
is constant, say $c$. Because, for every vector field $X$ on $M$, at each point in $M$, 
\begin{align}
X\,\vert\tau(\varphi)\vert^2=2\langle\overline{\nabla}_X\tau(\varphi),\tau(\varphi)\rangle=0.
\end{align}
\par
Therefore, if ${\rm Vol}(M,g)=\infty$ and $c\not=0$, 
then 
\begin{equation}
\tau_2(\varphi)=\frac12\int_M\vert \tau(\varphi)\vert^2\,v_g=\frac{c^2}{2}\,{\rm Vol}(M,g)=\infty
\end{equation}
which yields a contradiction. Thus, we have $\vert \tau(\varphi)\vert=c=0$, i.e., $\varphi$ is harmonic.  We have $(2)$. 
\vskip0.6cm\par
({\it The fourth step}) \quad 
For 
$(1)$, assume both 
$E(\varphi)<\infty$ and $E_2(\varphi)<\infty$. 
Then, let us consider a $1$-form 
$\alpha$ on $M$ defined by 
\begin{equation}
\alpha(X):=
\langle d\varphi(X),\tau(\varphi)\rangle, 
\qquad (X\in {\frak X}(M)). 
\end{equation}
Note here that 
\begin{align}
\int_M\vert\alpha\vert\,v_g
&=\int_M\left(
\sum_{i=1}^m\vert \alpha(e_i)\vert^2\right)^{1/2}\,v_g\nonumber\\
&\leq
\int_M\vert d \varphi\vert\,
\vert\tau(\varphi)\vert\,
v_g\nonumber\\
&\leq
\left(\int_M
\vert d\varphi\vert^2\,v_g\right)^{1/2}\,
\left(\int_M\vert\tau(\varphi)\vert^2\,v_g\right)^{1/2}
\nonumber\\
&=2\,\sqrt{E(\varphi)\,E_2(\varphi)}<\infty. 
\end{align}
Moreover, the divergent 
$\delta\alpha:=-\sum_{i=1}^m(\nabla_{e_i}\alpha)(e_i)\in C^{\infty}(M)$ turns out (cf. \cite{EL1}, 
p. 9) 
that 
\begin{equation}
-\delta\alpha=
\vert\tau(\varphi)\vert^2+\langle d\varphi,\overline{\nabla}\tau(\varphi)\rangle
=\vert\tau(\varphi)\vert^2.
\end{equation}
Indeed, we have 
\begin{align}
-\delta\alpha&=
\sum_{i=1}^me_i\langle d\varphi(e_i),\tau(\varphi)\rangle
-\sum_{i=1}^m\langle d\varphi(\nabla_{e_i}e_i),\tau(\varphi)\rangle\nonumber\\
&=\langle
\sum_{i=1}^m
\left(
\overline{\nabla}_{e_i}(d\varphi(e_i))-d\varphi(\nabla_{e_i}e_i)
\right),\tau(\varphi)\rangle\nonumber\\
&\qquad+\sum_{i=1}^m
\langle d\varphi(e_i),\overline{\nabla}_{e_i}\tau(\varphi)\rangle\nonumber\\
&=\langle\tau(\varphi),\tau(\varphi)\rangle
+\langle d\varphi,\overline{\nabla}\tau(\varphi)\rangle\nonumber
\end{align}
which is equal to 
$\vert\tau(\varphi)\vert$ 
since $\overline{\nabla}\tau(\varphi)=0$. 
\par
By $(3.16)$ and $E_2(\varphi)=\frac12\int_M\vert\tau(\varphi)\vert^2\,v_g<\infty$, 
the function $-\delta\alpha$ is also integrable over $M$. Thus, together with $(3.15)$, 
we can apply Gaffney's theorem (see 5.1 in Appendices, below) 
for the $1$-form $\alpha$.  
Then, by integrating $(3.16)$ over $M$,
and by Gaffney's theorem, we have 
\begin{equation}
0=\int_M(-\delta\alpha)\,v_g
=\int_M\vert\tau(\varphi)\vert^2\,v_g,
\end{equation}
which yields that $\tau(\varphi)=0$. 
\par
We have Theorem 2.1. 
\qed
\vskip0.6cm\par
Our method can be applied 
to an isometric immersion
 $\varphi:\,(M,g)\rightarrow (N,h)$. 
In this case,  
the $1$-form $\alpha$ defined by $(3.14)$ 
in the proof of Theorem 2.3 
vanishes automatically without using 
Gaffney's theorem 
since $\tau(\varphi)=m\,\xi$ belongs to the 
normal component 
of $T_{\varphi(x)}N$ $(x\in M)$, 
where $\xi$ is the mean curvature normal 
vector field and $m=\dim (M)$. 
Thus,  $(3.16)$ turns out that 
\begin{equation}
0=-\delta\alpha
=\vert\tau(\varphi)\vert^2+
\langle d\varphi,\overline{\nabla}\tau(\varphi)\rangle
=\vert\tau(\varphi)\vert^2
\end{equation}
which implies that $\tau(\varphi)=m\,\xi=0$, i.e., 
$\varphi$ is minimal. Thus, we obtain
\begin{th}
Assume that $(M,g)$ is a complete Riemannian manifold, and let 
$\varphi:\,(M,g)\rightarrow (N,h)$ 
is an isometric immersion, and the sectional curvature of $(N,h)$ is non-positive. 
If $\varphi:\,(M,g)\rightarrow (N,h)$ is biharmonic 
and $\int_M\vert\xi\vert^2\,v_g<\infty$, then 
$\varphi$ is minimal. 
Here, $\xi$ is the mean curvature normal vector field of the isometric immersion of $\varphi$. 
\end{th}
\vskip0.6cm\par
We also apply Theorem 2.1 to 
a horizontally conformal submersion 
$\varphi:\,(M^m,g)\rightarrow (N^n,h)$ $(m>n\geq 2)$ (cf. \cite{BW}, see also \cite{Gu}).   In the case that a Riemannian submersion from a space form of constant sectional curvature into a Riemann surface $(N^2,h)$, 
Wang and Ou (cf. \cite{WO}, see also \cite{LOu}) 
showed that it is 
biharmonic if and only if it is harmonic. We treat with 
a submersion from a higher dimensional Riemannian manifold $(M,g)$ (cf. \cite{BE}). 
Namely, let $\varphi: M\rightarrow N$ be a submersion, and 
each tangent space $T_xM$ $(x\in M)$ is decomposed into the orthogonal direct sum 
of the {\it vertical space} 
${\mathcal V}_x={\rm Ker}(d\varphi_x)$ and 
the {\it horizontal space} 
${\mathcal H}_x$: 
\begin{equation}
T_xM={\mathcal V}_x\oplus {\mathcal H}_x, 
\end{equation} 
and we assume that 
there exists a positive $C^{\infty}$ function 
$\lambda$ on $M$, 
called the {\it dilation}, such that, for each $x\in M$, 
\begin{equation}
h(d\varphi_x(X),d\varphi_x(Y))=\lambda^2(x)\,g(X,Y), \quad (X,Y\in {\mathcal H}_x). 
\end{equation}
The map $\varphi$ is said to be {\it horizontally homothetic} if the dilation $\lambda$ is constant along horizontally curves in $M$.  
\par
If
$\varphi:\,(M^m,g)\rightarrow (N^n,h)$ 
$(m>n\geq 2)$ is a horizontally conformal submersion . 
Then, the tension field $\tau(\varphi)$ is given (cf. \cite{BE}, \cite{BW}) by  
\begin{equation}
\tau(\varphi)=\frac{n-2}{2}\,\lambda^2\,d\varphi
\big({\rm grad}_{\mathcal H}\big(\frac{1}{\lambda^2}\big)\big)-(m-n)d\varphi\big(\hat{\bf H}\big), 
\end{equation}
where 
${\rm grad}_{\mathcal H}\big(\frac{1}{\lambda^2}\big)$ is the ${\mathcal H}$-component 
of the decomposition according to $(3.19)$ of 
${\rm grad}\big(\frac{1}{\lambda^2}\big)$, and 
$\hat{\mathbf {H}}$ is the trace of the second fundamental form of each fiber which is given by 
$\hat{\bf H}=
\frac{1}{m-n}\sum_{k=n+1}^m {\mathcal H}(\nabla_{e_k}e_k)$, 
where 
a local orthonormal frame field 
$\{e_i\}_{i=1}^m$ on $M$ 
is taken in such a way that 
$\{e_i{}_x\vert i=1,\cdots,n\}$ belong to 
${\mathcal H}_x$ and  
$\{e_j{}_x\vert j=n+1,\cdots,m\}$ belong to 
${\mathcal V}_x$ where $x$ is in a neighborhood in $M$. 
Then, due to Theorems 2.1 and (3.21), we have immediately
\begin{th}
Let $(M^m,g)$ be a complete non-compact Riemannian manifold, and 
$(N^n,h)$, a Riemannian manifold with the non-positive sectional curvature $(m>n\geq2)$.  
Let 
$\varphi:\,(M,g)\rightarrow (N,h)$ be a horizontally conformal submersion with the dilation $\lambda$  satisfying that 
\begin{equation}
\int_M\lambda^2\,\bigg\vert
\frac{n-2}{2}\,\lambda^2\,{\rm grad}_{\mathcal H}\big(\frac{1}{\lambda^2}\big)-(m-n)\,\hat{\bf H}\bigg\vert^{\,\,\,2}_g\,v_g<\infty.
\end{equation} 
Assume that, either
$\int_M\lambda^2\,v_g<\infty$ or 
${\rm Vol}(M,g)=\int_M\,v_g=\infty$. 
Then, if $\varphi:\,(M,g)\rightarrow (N,h)$ is biharmonic, then it is a harmonic morphism.  
\end{th}
\vskip0.6cm\par
Due to Theorem 3.2, we have: 
\begin{cor}
Let $(M^m,g)$ be a complete non-compact  Riemannian manifold, and 
$(N^2,h)$, a Riemannian surface with the non-positive sectional curvature $(m>n=2)$.  
Let 
$\varphi:\,(M,g)\rightarrow (N,h)$ be a horizontally conformal submersion with the dilation $\lambda$ satisfying that 
\begin{equation}
\int_M\lambda^2\,\big\vert
\hat{\bf H}\big\vert^{\,\,\,2}_g\,v_g<\infty.
\end{equation} 
Assume that, either
$\int_M\lambda^2\,v_g<\infty$ or 
${\rm Vol}(M,g)=\int_M\,v_g=\infty$. 
Then, if $\varphi:\,(M,g)\rightarrow (N,h)$ is biharmonic, then it is a harmonic morphism.  
\end{cor}
\vskip0.6cm\par
Corollary 3.3 implies 
\begin{cor}
Let $(M^m,g)$ be a non-compact complete Riemannian manifold $(m>2)$, and $(N^2,h)$, a Riemannian surface with non-positive curvature.  
Let  $\lambda$ be a positive function  
in $C^{\infty}(M)\,\cap \,L^2(M)$, 
where $L^2(M)$ is the space of square integrable functions on $(M,g)$. 
Then,  
every 
biharmonic 
horizontally conformal submersion 
$\varphi:\,(M^m,g)\rightarrow (N^2,h)$ 
with  a dilation $\lambda$ and a 
bounded $\vert\hat{\bf H}\vert_g$,  
exactly 
$\lambda\,\vert\hat{\bf H}\vert_g\in L^2(M)$, 
must be 
a harmonic morphism. 
\end{cor}
\vskip0.6cm\par
\begin{rem} 
Notice that in Corollary 3.4, (1), there is no restriction to the dilation $\lambda$ because of 
$\dim N=2$. This implies that for every 
positive $C^{\infty}$ function $\lambda$ in 
$C^{\infty}(M)\,\cap\,L^2(M)$ 
satisfying 
(3.2), 
we have 
a harmonic morphism 
$\varphi: \,(M^m,g)\rightarrow (N^2,h)$. 
\end{rem}
\vskip2cm\par
\section{Appendices}
\subsection{Gaffney's theorem}
In this appendix, we recall Gaffney's theorem (\cite{G}):
\begin{th} (Gaffney)
Let $(M,g)$ be a complete Riemannian manifold. 
If a $C^1$ $1$-form $\alpha$ satifies that 
$\int_M\vert\alpha\vert\,v_g<\infty$ and 
$\int_M(\delta\alpha)\,v_g<\infty$, 
or equivalently, 
a $C^1$ vector field $X$ defined by 
$\alpha(Y)=\langle X,Y\rangle$ $(\forall\,\,Y\in {\frak X}(M))$ 
satisfies that 
$\int_M\vert X\vert \,v_g
<\infty$
and $\int_M\text{\rm div}(X)\,v_g<\infty$, 
then
\begin{equation}
\int_M(-\delta\alpha)\,v_g=\int_M{\text{\rm div}}(X)\,v_g=0.
\end{equation}
\end{th}
\begin{pf} 
For completeness, we give a proof. 
By integrating over $M$, the both hand sides of 
\begin{equation}
\text{\rm div}(\eta^2\,X)=\eta^2\,\text{\rm div}(X)
+2\eta\,\langle\nabla\eta,X\rangle, 
\end{equation}
we have 
\begin{equation}
\int_M\text{div}(\eta^2\,X)\,v_g
=\int_M\eta^2\,\text{\rm div}(X)\,v_g
+2\int_M\eta\,\langle\nabla\eta,X\rangle\,v_g. 
\end{equation}
Since the support of $\eta^2\,X$ is compact, 
the left hand side must vanish. So, we have
\begin{equation}
\int_M
\eta^2\,\text{\rm div}(X)\,v_g
=-2\int_M\eta\,\langle\nabla\eta,X\rangle\,v_g. 
\end{equation}
Therefore, we have
\begin{align}
\bigg\vert\int_{B_r(x_0)}\text{div}(X)\,v_g\bigg\vert
&\leq
\bigg\vert\int_M\eta^2\,\text{\rm div}(X)\,v_g
\bigg\vert\nonumber\\
&=2\,\bigg\vert
\int_M\eta\,\langle\nabla\eta,X\rangle\,v_g\bigg\vert\nonumber\\
&
\leq 2\int_M\eta\,\vert\nabla\eta\vert\,\vert X\vert\,v_g\nonumber\\
&\leq
\frac{4}{r}\int_M\vert X\vert\,v_g. 
\end{align}
By the assumption that $\int_M\vert X\vert\,v_g<\infty$, 
the right hand side goes to $0$ if $r$ tends to infinity.  
Since $B_r(x_0)$ goes to $M$ 
as $r\rightarrow\infty$, due to completeness of 
$(M,g)$, and the assumption 
that $\int_M\text{\rm div}(X)\,v_g<\infty$, 
we have
$\int_M\text{div}(X)\,v_g=\lim_{r\rightarrow\infty}
\int_{B_r(x_0)}\text{\rm div}(X)\,v_g=0$. 
\end{pf}
\vskip0.6cm\par
\subsection{The Bochner-type estimations and alternative proof of main theorem}
We first show the Bochner-type estimations for the tension fields 
of biharmonic maps into a Riemannian manifold $(N,h)$ of 
non-positive curvature. 
We write down them for clarity since 
at least we do not find them in the literature. 
Then we can give an alternative proof of our main theorem (cf. Theorem 2.1) by using them.  
\begin{lem}
Assume that the sectional curvature of 
$(N,h)$ is non-positive, 
and $\varphi:\,(M,g)\rightarrow (N,h)$ is a biharmonic mapping. 
Then, it holds that 
\begin{equation}
\Delta\,\vert\tau(\varphi)\vert^2
\geq 
2\,\vert\overline{\nabla}\tau(\varphi)\vert^2
\end{equation}
in $M$. Here, $\Delta=\sum_{i=1}^m(e_i{}^2-\nabla_{e_i}e_i)$ 
is 
the Laplace-Beltrami operator of 
$(M,g)$. 
\end{lem}
\begin{pf}
\quad Let us take a local orthonormal frame field 
$\{e_i\}_{i=1}^m$ 
on $M$, and 
$\varphi:\,(M,g)\rightarrow (N,h)$, a biharmonic map. 
Then, for $V:=\tau(\varphi)\in \Gamma(\varphi^{-1}TN)$, we have 
\begin{align}
\frac12\,\Delta\,\vert V\vert^2
&=\frac12
\sum_{i=1}^m
\left\{
e_i{}^2\,\vert V\vert^2-\nabla_{e_i}e_i\,\vert V\vert^2
\right\}\nonumber\\
&=
\sum_{i=1}^m
\left\{
e_i\,h(\overline{\nabla}_{e_i}V,V)
-h(\overline{\nabla}_{\nabla_{e_i}e_i}V,V)
\right\}\nonumber\\
&=
\sum_{i=1}^m
\left\{
h(\overline{\nabla}_{e_i}\overline{\nabla}_{e_i}\,V,V)
-h(\overline{\nabla}_{\nabla_{e_i}e_i}\,V,V)
\right\}
\nonumber\\
&\qquad
+\sum_{i=1}^mh(\overline{\nabla}_{e_i}V,
\overline{\nabla}_{e_i}V)
\nonumber\\
&=h(-\overline{\Delta}\,V,V)
+\vert\overline{\nabla}V\vert^2\nonumber\\
&=h(-{\mathcal R}(V),V)+\vert\overline{\nabla}V\vert^2\nonumber\\
&\geq\vert\overline{\nabla}V\vert^2,
\end{align}
because for the second last equality, we used $\overline{\Delta}V-{\mathcal R}(V)=J(V)=0$ for $V=\tau(\varphi)$, 
due to the biharmonicity of $\varphi:\,(M,g)\rightarrow (N,h)$, 
and for the last inequality of $(4.7)$, 
we used
\begin{equation}
h({\mathcal R}(V),V)
=\sum_{i=1}^m
h(R^N(V,\varphi_{\ast}e_i)\varphi_{\ast}e_i,V)\leq 0
\end{equation}
since the sectional curvature of $(N,h)$ is non-positive. 
\end{pf}
\vskip0.6cm\par
By Lemma 4.2, we have 
\begin{lem}
Under the same assumptions as Lemma 4.2, we have 
\begin{equation}
\vert\tau(\varphi)\vert\,\,\Delta\,\vert\tau(\varphi)\vert\geq 0.
\end{equation}
\end{lem}
\begin{pf}
Due to Lemma 4.2, we have 
\begin{align}
2\,\vert\overline{\nabla}\tau(\varphi)\vert^2
&\leq \Delta\,\vert\tau(\varphi)\vert^2\nonumber\\
&=2\,\vert\tau(\varphi)\vert\,\,\Delta\,\vert\tau(\varphi)\vert
+2\,\vert\,\nabla\,\vert\tau(\varphi)\vert\,\,\vert^2.
\end{align}
Thus, we have 
\begin{align}
\vert\tau(\varphi)\vert\,\Delta\,\vert(\tau(\varphi)\vert
&\geq
\vert\overline{\nabla}\tau(\varphi)\vert^2-\vert\,\nabla\,\vert\tau(\varphi)\vert\,\,\vert^2\nonumber\\
&\geq 0 
\end{align}
on the set of points in $M$ 
where $\vert\tau(\varphi)\vert>0$. 
Here, to see the last inequality of (4.11), it suffices to 
notice that for all $V\in \Gamma(\varphi^{-1}TN)$, 
\begin{equation}
\vert\overline{\nabla}V\vert\geq \vert\,\,\nabla\,\vert V\vert\,\,\vert 
\end{equation}
on the set of points in $M$ 
where $\vert V\vert > 0$, 
which follows from that 
\begin{align}
\vert V\vert\,\,\vert\,\nabla\,\vert V\vert\,\,\vert
&=\frac12\,\vert\,\nabla\,\vert V\,\vert^2\,\,\vert\nonumber\\
&=\frac12\,\vert \nabla\,h(V,V)\,\vert\nonumber\\
&=\vert\,h(\overline{\nabla}V,V)\,\vert\nonumber\\
&\leq \vert\,\overline{\nabla}V\,\vert\,\,\vert V\vert.
\end{align}
This proves Lemma 4.3. 
\end{pf}
\vskip0.6cm\par
Now,  we are in position to give an alternative  proof of  Theorem 2.1. 
\vskip0.6cm\par
({\it The first step}) \quad 
For a fixed point $x_0\in M$, and 
for every 
$0<r<\infty$, 
we first take the same  cutoff  $C^{\infty}$ function $\eta$ on $M$ in Section Three. 
\par
For $0<r<\infty$, multiply 
$\lambda^2$ to both hand sides of the inequality (4.9) in Lemma 4.3,
and integrate over $M$, we have 
\begin{align}
0&\leq\int_M
\eta^2\,\vert\tau(\varphi)\vert\,\Delta\,(\vert\tau(\varphi)\vert)\, v_g
\nonumber\\
&=-\int_M\langle\,\nabla(\eta^2\,\vert\tau(\varphi)\vert
),\nabla\,\vert\tau(\varphi)\vert\,\rangle\,v_g\nonumber\\
&=
-\int_M
\eta^2\,\vert\,\nabla(\vert\tau(\varphi)\vert)\,\vert^2\,v_g\nonumber\\
&\quad -2\int_M\vert\tau(\varphi)\vert\,\eta\,
\langle\,\nabla(\vert\tau(\varphi)\vert),\nabla\eta\,\rangle\,v_g\nonumber\\
&=
-\int_M\vert \nabla(\vert\tau(\varphi)\vert)\vert^2\,\eta^2\,v_g\nonumber\\
&\quad-2\int_M
\langle\,\eta\,\nabla(\vert\tau(\varphi)\vert),
\vert\tau(\varphi)\vert\,\nabla\eta\,\rangle
\,v_g.
\end{align} 
Therefore, by applying 
$A:=\eta\,\nabla(\vert\tau(\varphi)\vert)$ and 
$B:=\vert\tau(\varphi)\vert\,\nabla\eta$ 
to Young's inequality:\,
$-2A\,B\leq\epsilon^2A^2+\frac{1}{\epsilon^2}B^2$, 
for every positive real number 
$\epsilon>0$, 
the right hand side of $(4.14)$ is smaller than or equal to 
\begin{align}
&-\int_M\vert \nabla\,(\vert\tau(\varphi)\vert)\vert^2\,\eta^2\,v_g\nonumber\\
&+\epsilon^2\int_M\vert\nabla\,(\vert\tau(\varphi)\vert)\,\vert^2\,\eta^2\,v_g+
\frac{1}{\epsilon^2}\int_M\vert\tau(\varphi)\vert^2\,
\vert\nabla\eta\vert^2\,v_g\nonumber\\
&=-(1-\epsilon^2)
\int_M\vert\nabla\,(\vert\tau(\varphi)\vert)\,\vert^2\,\eta^2\,v_g+
\frac{1}{\epsilon^2}\int_M\vert\tau(\varphi)\vert^2\,
\vert\nabla\eta\vert^2\,v_g.
\end{align}
Thus, we obtain 
\begin{equation}
(1-\epsilon^2)
\int_M\vert\nabla\,(\vert\tau(\varphi)\vert)\,\vert^2\,\eta^2\,v_g
\leq
\frac{1}{\epsilon^2}\int_M\vert\tau(\varphi)\vert^2\,
\vert\nabla\eta\vert^2\,v_g.
\end{equation}
\vskip0.6cm\par
({\it The second step}) \quad 
In the inequality $(4.16)$, the left hand side 
is bigger than or equal to 
$(1-\epsilon^2)\int_{B_r(x_0)}
\vert\nabla\,(\vert\tau(\varphi)\vert)\,\vert^2\,\eta^2\,v_g$ 
since $\eta=1$ on $B_r(x_0)$, 
and 
the right hand side is smaller than or equal to 
$
\frac{1}{\epsilon^2}\int_M\vert\tau(\varphi)\vert^2\frac{4}{r^2}\,v_g
=\frac{4}{\epsilon^2\,r^2}\int_M\vert\tau(\varphi)\vert^2\,v_g
$, 
we obtain 
\begin{equation}
(1-\epsilon^2)\int_{B_r(x_0)}
\vert\nabla\,(\vert\tau(\varphi)\vert)\,\vert^2\,\eta^2\,v_g
\leq
\frac{4}{\epsilon^2\,r^2}\int_M\vert\tau(\varphi)\vert^2\,v_g. 
\end{equation}
By putting $\epsilon=\frac12$, we have 
\begin{equation}
\int_{B_r(x_0)}
\vert\nabla\,(\vert\tau(\varphi)\vert)\,\vert^2\,\eta^2\,v_g
\leq
\frac{64}{3\,r^2}\int_M\vert\tau(\varphi)\vert^2\,v_g. 
\end{equation}
\vskip0.6cm\par
({\it The third step}) \quad 
Now we may take $r\rightarrow \infty$ 
in $(4.18)$, 
due to the assumptions that 
$\int_M\vert\tau(\varphi)\vert^2\,v_g<\infty$, 
and $(M,g)$ is complete, 
we obtain 
\begin{equation}
\int_M
\vert\nabla\,(\vert\tau(\varphi)\vert)\,\vert^2\,v_g
=0,
\end{equation}
which implies $\nabla\,(\vert\tau(\varphi)\vert)=0$, 
that is, 
\begin{equation}
\vert\tau(\varphi)\vert=c \quad (\text{\rm a constant})
\end{equation} 
on $M$. 
Then, the inequality $(4.6)$ turns out that
\begin{equation}
0=
\Delta(\vert\tau(\varphi)\vert^2)
\geq
2\,\vert\overline{\nabla}\tau(\varphi)\vert^2
\geq
0,
\end{equation} 
which means that 
\begin{equation}
\overline{\nabla}\tau(\varphi)=0 \qquad(\text{on}\,\,M).
\end{equation}
\vskip0.6cm\par
By passing through the same procedure of 
the fourth step of the proof of Theorem 2.1 in Section Three,  we have done. 
\qed
\newpage\par

\end{document}